\documentclass[12pt]{article}

\usepackage{amssymb}
\usepackage{amsmath}
\usepackage{amsfonts}
\usepackage{setspace}
\usepackage{epsfig}
\usepackage{subfigure}
\usepackage{color}
\usepackage{enumerate}

\numberwithin{equation}{section}

\begin{document}
\baselineskip=22pt

\newcommand{\ve}[1]{\mbox{\boldmath$#1$}}
\newcommand{\be}{\begin{equation}}
\newcommand{\ee}{\end{equation}}
\newcommand{\bc}{\begin{center}}
\newcommand{\ec}{\end{center}}
\newcommand{\bal}{\begin{align*}}
\newcommand{\enal}{\end{align*}}
\newcommand{\al}{\alpha}
\newcommand{\bt}{\beta}
\newcommand{\gm}{\gamma}
\newcommand{\de}{\delta}
\newcommand{\la}{\lambda}
\newcommand{\om}{\omega}
\newcommand{\Om}{\Omega}
\newcommand{\Gm}{\Gamma}
\newcommand{\De}{\Delta}
\newcommand{\Th}{\Theta}
\newcommand{\kp}{\kappa}
\newcommand{\nno}{\nonumber}
\newtheorem{theorem}{Theorem}[section]
\newtheorem{lemma}{Lemma}[section]
\newtheorem{assum}{Assumption}[section]
\newtheorem{claim}{Claim}[section]
\newtheorem{proposition}{Proposition}[section]
\newtheorem{corollary}{Corollary}[section]
\newtheorem{definition}{Definition}[section]
\newtheorem{remark}{Remark}[section]
\newenvironment{proof}[1][Proof]{\begin{trivlist}
\item[\hskip \labelsep {\bfseries #1}]}{\end{trivlist}}
\newenvironment{proofclaim}[1][Proof of Claim]{\begin{trivlist}
\item[\hskip \labelsep {\bfseries #1}]}{\end{trivlist}}

\def \qed {\hfill \vrule height7pt width 5pt depth 0pt}
\def\refhg{\hangindent=20pt\hangafter=2pt}
\def\refmark{\par\vskip 2.50mm\noindent\refhg}
\def\Bbb#1{\mbox{\sf #1}}

\title{\textbf{Super Fractal Interpolation Functions}}
\date{}
\author{G.P. Kapoor$^1$  and Srijanani Anurag Prasad$^2$}
\maketitle
\vspace{-1.8cm}
\bc Department of Mathematics and Statistics\\
Indian Institute of Technology Kanpur\\  Kanpur 208016 India\\
$^1$gp@iitk.ac.in    $^2$jana@iitk.ac.in\ec

\doublespacing

\begin{abstract}
In the present work, the notion of Super Fractal Interpolation
Function (SFIF) is introduced for finer simulation of the objects of
the nature or outcomes of scientific experiments that reveal one or
more structures embedded in to another.  In the construction of
SFIF,  an IFS is chosen from a pool of several IFS at each level of
iteration leading to implementation of the desired randomness and
variability in fractal interpolation of the given data. Further, an
expository description of our investigations on the integral, the
smoothness and determination of conditions for existence of
derivatives  of a SFIF is given in the present work.
\end{abstract}



\newpage

\section{Introduction}

Barnsley~\cite{barnsley86}
introduced Fractal Interpolation Function (FIF) using the theory of Iterated Function System (IFS). Since then, a growing  numbers of papers have been published
showing relation between fractals and wavelets~\cite{donovan96,guo10}, fractal functions and Kiesswetter-like functions~\cite{yong04} and on fractal dimension~\cite{yanyan03,liang07}. Later, Barnsley et. al. extended the idea of FIF to
produce more flexible interpolation functions called Hidden-variable FIF (HFIF) which were generally non-self affine. Dalla~\cite{dalla03} found bounds on fractal dimension for the graphs of non-affine FIFs. In 1989, Barnsley
and Harrington~\cite{barnsley89_1} constructed an IFS to show that a FIF can be indefinitely integrated, giving rise to a
hierarchy of smoother functions and developed results on differentiability of a FIF.

Fractal Interpolation Function, constructed as attractor
of a single Iterated Function System (IFS) by virtue of
self-similarity alone, is not rich enough to describe an object
found in nature or output of a certain scientific experiment.
The objects of nature generally reveal one or more structures embedded
in to another. Similarly, the outcomes of several scientific
experiments exhibit randomness and variation at various stages.
Therefore, more than one IFSs are needed to model such objects.
Barnsley~\cite{barnsley05,barnsley06,barnsley08} introduced the class of super
fractal sets constructed by using multiple IFSs to simulate such
objects. Massopust~\cite{massopust10} constructed super fractal functions and V-variable fractal functions by joining pieces of fractal functions which are attractor of finite family of IFss.  However, for a data set arising from nature or a scientific
experiment, a solution of fractal interpolation problem based on
several IFS has not been investigated so far. To fill this gap, the
notion of Super Fractal Interpolation Function~(SFIF) is introduced
in the present work. The construction of SFIF requires the use of
more than one IFS wherein, at each level of iteration, an IFS can be
chosen from a pool of several IFS. This approach is likely to ensure
desired randomness and variability needed to facilitate better
geometrical modeling of objects found in nature and results of
certain scientific experiments. The construction of SFIF is followed
in the present paper by the investigations of its smoothness, its
integral  and determination of conditions for existence of its
derivatives.

The organization of the present paper is as follows: In
Section~\ref{sec:sfif},  for a given finite set of data, the method
of construction of a Super Fractal Interpolation Function (SFIF) is
developed. At each level of iteration, an IFS is chosen from a pool
of IFS in our construction of SFIF.  For a sample interpolation
data, a computational model of SFIF, illustrating the construction
method given in Section~\ref{sec:sfif},  is presented in
Section~\ref{sec:exsfif}. The fractal dimension and average fractal
distance  are computed for various SFIFs constructed in this
section. Finally, in Section~\ref{sec:idsfif}, it is found that for
a SFIF passing through a given interpolation data, its integral is
also a SFIF, albeit for a different interpolation data. An
expository description of  smoothness of a SFIF and conditions for
existence of derivatives of a SFIF is also given in this section.

\section{Construction of SFIF}\label{sec:sfif}

In this section,  the  notion of Super Fractal Interpolation
Function (SFIF) is introduced via its construction based on more
than one IFS.

Let $ S_0 = \{ (x_i,y_i) \in \mathbb{R}^2 :i=0, \ldots N \}$ be the
set of given interpolation data.  For $k=1,\ldots,M$, $M>1$ and $ n
\ = \ 1,\ldots,N $, let the functions $\ \om_{n,k} :I \times
\mathbb{R} \rightarrow I \times \mathbb{R} $ be defined by
\begin{align}\label{eq:wnk}
 \om_{n,k}(x,y)=(L_n(x),G_{n,k}(x,y)) \
\mbox{for all} \ (x,y) \in \mathbb{R}^2
\end{align}
where, the contractive homeomorphisms $ L_n : I \rightarrow I_n $
are given by
\begin{align}\label{eq:Ln}
L_n(x) = a_n x + b_n = \frac{(x_n-x_{n-1})  x  +(x_N x_{n-1}-x_0
x_n)}{(x_N-x_0)}
\end{align}
and the functions $G_{n,k} : I \times \mathbb{R} \rightarrow
\mathbb{R}$ defined by
\begin{align}\label{eq:gnk}
G_{n,k}(x,y)= e_{n,k} x + \gm_{n,k} y +  f_{n,k}
\end{align}
satisfy the join-up conditions
\begin{align}\label{eq:jc}
G_{n,k}(x_0, y_{0})& =  y_{n-1} \quad \mbox{and} \quad G_{n,k}(x_N,
y_{N}) = y_n.
\end{align}
Here, $\gm_{n,k}$ are free parameters chosen such that $|\gm_{n,k}|
< 1$ and $\gm_{n,k} \neq \gm_{n,l}$ for $k \neq l$.
By~\eqref{eq:jc}, it is observed that $\ \om_{n,k}$ are continuous
functions. The Super Iterated Function System (SIFS) that is needed
to construct SFIF corresponding to the set of given interpolation
data $S_0 = \{ (x_i,y_i) \in \mathbb{R}^2 :i=0,1 \ldots,N  \}$ is
now defined as the pool of IFS
\begin{align}\label{eq:sifs}
\Big\{ \big\{\mathbb{R}^2; \om_{n,k} : n =1, \ldots,N \big\}, \ k=
1,\ldots,M \Big\}
\end{align}
where, the functions $\om_{n,k}$ are given by~\eqref{eq:wnk}.

To introduce a SFIF associated with SIFS~\eqref{eq:sifs},  let
$\{W_k : \mathcal{H}(\mathbb{R}^2) \rightarrow
\mathcal{H}(\mathbb{R}^2),   k=1,\ldots,M \} $,  be a collection of
continuous functions defined by $ W_k(G) = \bigcup\limits_{n=1}^N
\om_{n,k}(G)\  \mbox{where}, \ \om_{n,k}(G) = \om_{n,k}(x,y) \
\mbox{for all} \ (x,y) \in G$.  Since, $ h(W_k(A),W_k(B)) \leq
\max\limits_{1 \leq n \leq N} \gm_{n,k} \ h(A,B)$, where $h$ is
Hausdorff metric on $\mathcal{H}(\mathbb{R}^2)$, $
\left\{\mathcal{H}(\mathbb{R}^2);\ W_1,\ldots, W_M \right\} $ is a
hyperbolic IFS. Hence, by Banach fixed point theorem, there exists
an attractor $\mathcal{A} \in
\mathcal{H}(\mathcal{H}(\mathbb{R}^2))$.

 Let $ \Lambda $ be the code space on $M$ natural
numbers $ {1,2,\ldots,M}$ . For $\sigma = \sigma_1\sigma_2\ldots
\sigma_k\ldots\in \Lambda$, define the function $\phi : \Lambda
\rightarrow \mathcal{H}(\mathbb{R}^2)$  by
\begin{align}\label{eq:phi}
 \phi(\sigma) = \lim\limits_{k \rightarrow \infty} W_{\sigma_k}
\circ W_{\sigma_{k-1}} \circ \ldots \circ W_{\sigma_1}(G),\ G \in
\mathcal{H}(\mathbb{R}^2) .
\end{align}
It is shown that~\cite{barnsley88} $\phi(\sigma)$ exists, belongs to
$\mathcal{A}$ and is independent of $G \in
\mathcal{H}(\mathbb{R}^2)$.  Also,  the function $\phi$ is onto and
continuous~\cite{barnsley88}. In the construction of SFIF, for a
$\sigma =\sigma_1\sigma_2 \ldots \in \Lambda$, let  the action of
SIFS~\eqref{eq:sifs} at the iteration level $j$ be defined by $S_j =
W_{\sigma_j}(S_{j-1})$, where $S_0 $ is the set of given
interpolation data. It is easily seen that the set,
\begin{align}\label{eq:Gsigma}
G_{\sigma} \equiv \phi(\sigma) =~\lim\limits_{k \rightarrow \infty}
W_{\sigma_k} \circ \ldots \circ W_{\sigma_1} (S_0) = \lim\limits_{k
\rightarrow \infty} S_k
\end{align}  is
the attractor of SIFS~\eqref{eq:sifs} for a fixed $\sigma \in
\Lambda$. The following theorem shows that $G_{\sigma}$ is the graph
of a continuous function~$g_{\sigma}$.

\begin{theorem}\label{th:SFIF}
Let $G_{\sigma} $ be the attractor of SIFS~\eqref{eq:sifs} for
$\sigma = \sigma_1\sigma_2\ldots \sigma_k\ldots\in \Lambda$. Then,
$G_{\sigma}$ is graph of a continuous function $g_{\sigma} : I
\rightarrow \mathbb{R}$ such that $ g_{\sigma} (x_n) = y_n$ for all
$\ n=0,\ldots,N$.
\end{theorem}

\begin{proof}
Let $g_0$ be a function whose graph is $S_0$.  Then, the set $S_k$,
$k \geq 1$, is graph of the function $ g_{\sigma_k}$, where $
g_{\sigma_k} (x) = G_{i_k,\sigma_{k}} \big( L_{i_k}^{-1}
(x),g_{\sigma_{k-1}}(L_{i_k}^{-1}(x))\big) $. It is easily seen that
$ g_{\sigma_k} (x) = G_{i_k,\sigma_{k}} \bigg( L_{i_k}^{-1} (x),
G_{i_{k-1},\sigma_{k-1}}\Big( . ,\ldots
G_{i_1,\sigma_1}(L_{i_1}^{-1}\circ \ldots \circ
L_{i_k}^{-1}(x),g_0(L_{i_1}^{-1}\circ \ldots \circ L_{i_k}^{-1}(x)))
\ldots \Big)\bigg) $ . Therefore, it follows by~\eqref{eq:Gsigma}
that the set $G_{\sigma} $ is graph of the function $ g_{\sigma} =
\lim\limits_{k \rightarrow \infty} g_{\sigma_k} $.

For proving the continuity of the function $ g_{\sigma}$, consider
$\tau_1^*\tau_2^* \ldots\tau^*_j\ldots \in \Lambda$ where $\tau^*_j
\neq 1$ for some $j \in \mathbb{N}$ and $\tau^*_i = 1$ for $ i \in
\mathbb{N}$ and $ i \neq j$. We first show that $G_{\tau^*} $ is
graph of a continuous function $g_{\tau^*}$. If not, then
$G_{\tau^*} = \phi(\tau^*) $ is graph of a function $g_{\tau^*}$
that is not continuous so that there exist a $\delta_1
> 0$ such that whenever $x_1, x_2 \in I$ and $|x_1 -x_2| < \delta_1$,
\begin{align} \label{eq:discn}
|g_{\tau^*}(x_1) - g_{\tau^*}(x_2)| > \epsilon.
\end{align}
It is known that~\cite{barnsley86}, for $\tau = \bar{1} \in
\Lambda$, $G_{\tau} = \phi(\tau) $, with $\phi$ defined
by~\eqref{eq:phi}, is graph of a continuous function $g_{\tau} : I
\rightarrow \mathbb{R}$ such that $ g_{\tau} (x_n) = y_n$, $
n=0,1,\ldots,N$. Consequently, there exists a $\delta_2 > 0$ such
that $|x_1 -x_2| < \delta_2$ implies $|g_{\tau}(x_1) -
g_{\tau}(x_2)| < \frac{\epsilon}{3}$. Also, since $\phi$ is a
continuous map, there exists $\delta_3 > 0$ such that, for $\tau$
and $\tau^*$ satisfying
 $d_c(\tau,\tau^*) = \frac{|\tau_j-\tau_j^*|}{(M+1)^j} < \delta_3$,
$\max\limits_{x \in I} |g_{\tau}(x) - g_{\tau^*}(x)| <
\frac{\epsilon}{3}$. Thus, for $\delta = \min (\de_1,\de_2,\de_3)$
and $x_1, x_2 $ satisfying $|x_1 -x_2| < \delta$,
 $|g_{\tau^*}(x_1) - g_{\tau^*}(x_2)| \leq |g_{\tau^*}(x_1) -
g_{\tau}(x_1)|+ |g_{\tau}(x_1) - g_{\tau}(x_2)| + |g_{\tau}(x_2) -
g_{\tau^*}(x_2)| <  \epsilon$,  a contradiction to~\eqref{eq:discn}.
Hence, $G_{\tau^*} $ is graph of continuous function~$g_{\tau^*}$.

 Now, consider the sequence
$\sigma_n=\sigma_{1,n}\sigma_{2,n}\ldots....$ with $\sigma_{j,n} =
\sigma_j $ for $j \leq n$ and $\sigma_{j,n} = 1$ for $ j > n$. It is
easily seen  that as $n$ tends to infinity, $\sigma_n$ tends to
$\sigma$ with respect to the metric $d_c$. Using the arguments of
previous paragraph inductively, it follows that $G_{\sigma_n} =
\phi(\sigma_n) $ is graph of a continuous function $g_{\sigma_n}$
defined on $I$. Let $G_{\sigma} = \phi(\sigma) $ be graph of a
function $g_{\sigma}$. By continuity of $\phi$, $G_{\sigma_n}$ tends
to $G_{\sigma}$ with respect to Hausdorff metric $h$ as $n
\rightarrow \infty$, which implies that $g_{\sigma_n}$ tends to
$g_{\sigma}$ with respect to Maximum metric  as $n \rightarrow
\infty$. Hence, there exist an $\epsilon
> 0$ such that $\max\limits_{x \in I} |g_{\sigma_n}(x) -
g_{\sigma}(x)| < \frac{\epsilon}{3} $. Since $g_{\sigma_n}$ is
continuous on $I$, there exits a $\de
> 0 $ such that $|x-y| < \de$ implies $|g_{\sigma_n}(x) -
g_{\sigma_n}(y)|~<~\frac{\epsilon}{3}$. Therefore,
$|g_{\sigma}(x)-g_{\sigma}(y)| \leq |g_{\sigma}(x)-g_{\sigma_n}(x)|
+ |g_{\sigma_n}(x) - g_{\sigma_n}(y)| + |g_{\sigma_n}(y)
-g_{\sigma}(y)| < \epsilon$ for $|x-y| < \de$ implying that the
function $g_{\sigma}$ is continuous on $I$. This establishes that
the attractor $G_{\sigma} $ of SIFS~\eqref{eq:sifs} is the graph of
continuous function~$g_{\sigma}$. \qed \end{proof}

Theorem~\ref{th:SFIF} is instrumental in defining a SFIF associated
with SIFS~\eqref{eq:sifs} as follows:

\begin{definition}
The  \textbf{Super Fractal Interpolation Function (SFIF) } for the
given interpolation data $\{(x_i,y_i) :i = 0,1,\ldots,N \}$   is
defined as the continuous function $g_{\sigma}$ whose graph
$G_{\sigma}$ is the attractor of SIFS~\eqref{eq:sifs}.
\end{definition}

\begin{remark}
Consider the family of continuous functions
$
 {\cal
G} = \{ f : I \rightarrow \mathbb{R} \ \mbox{such that} \\ f  \mbox{
is continuous}, f(x_0) = y_0 \ \mbox{and} \ f(x_N) = y_N\} $ with
metric  $ d_{{\cal G}} (f,g) = \max\limits_{x \in I } |f(x) - g(x)|
$. Since ${\cal G}$ is a complete metric space,  it is easily seen
that, for a fixed $\sigma \in \Lambda$,  Read-Bajraktarevic operator
$T : \Lambda \times {\cal G} \rightarrow {\cal G}$ defined as
\begin{align}\label{eq:rb}
T(\sigma, g)(x) & =  \lim\limits_{k \rightarrow \infty}
\bigg\{G_{i_k,\sigma_{k}} \bigg( L_{i_k}^{-1} (x),
G_{i_{k-1},\sigma_{k-1}} \Big( L_{i_{k-1}}^{-1}\circ L_{i_k}^{-1}
(x), G_{i_{k-2},\sigma_{k-2}}\big(. ,\ldots \nno \\ & \quad \quad
 \quad \quad G_{i_1,\sigma_1}(L_{i_1}^{-1}\circ \ldots \circ
L_{i_k}^{-1}(x),g(L_{i_1}^{-1}\circ \ldots \circ L_{i_k}^{-1}(x)))
\ldots \big)\Big)\bigg)\bigg\},
\end{align}
is a contraction map on ${\cal G}$ and so it has a unique fixed
point in ${\cal G}$. It is observed that, SFIF $g_{\sigma}$
satisfies $g_{\sigma} = T(\sigma, g_{\sigma})$.
\end{remark}

\begin{remark}\label{rem:kpsfif}
The notion of SFIF  can further be generalized by introducing a
fixed parameter  $\kappa \ (0 \leq \kappa <1)$ in the join-up
conditions~\eqref{eq:jc} as follows:
\begin{align*}\left. \begin{array}{ll}
G_{n,k}(x_0,\kp x_{0} + (1-\kp) y_{0})&=\kp x_{n-1} + (1-\kp)
y_{n-1} \\ G_{n,k}(x_N,\kp x_{N} + (1-\kp) y_{N})&=  \kp x_n + ( 1 -
\kp)y_n. \end{array} \right\}
\end{align*}
The above condition ensures that there exits a unique attractor
$G_{\sigma,\kp} \in {\cal H}(\mathbb{R}^2)$ of SIFS~\eqref{eq:sifs}.
By the arguments similar to those in the proof of
Theorem~\ref{th:SFIF}, $G_{\sigma,\kp}$ is graph of a continuous
 function $g_{\sigma,\kp}$, called henceforth  Parameterized
 SFIF or  $\kp$-SFIF.
\end{remark}

\section{Computational Model of SFIF}\label{sec:exsfif}

Our method of construction developed in Section~\ref{sec:sfif} is
employed in the present section
 for generating various SFIF for a sample interpolation data $S_0 = \{(0,0), (30,90 ),
(60,70),(100,10)\}$. For identifying the corresponding SIFS $ \Big\{
\big\{\mathbb{R}^2; \om_{n,k} : n =1, 2,3 \big\}, \ k= 1,2 \Big\} $,
the maps $\om_{n,k} ,\ k=1,2$ (c.f.~\eqref{eq:wnk}) are obtained by
computing (c.f. Table~\ref{tab:Table-coef}) the values of $a_i$,
$b_i$~(c.f.~\eqref{eq:Ln}) and $e_{i,1},\ f_{i,1};\ e_{i,2}$ and
$f_{i,2}$ (c.f.~\eqref{eq:jc}) with $\gm_{i,1} =0.4$ and $\gm_{i,2}
=0.6$ for $i=1,2,3$.

\begin{figure}[!htp]
{\centering \subfigure[ $\sigma_{(r)} =
\overline{112}$]{\epsfig{file=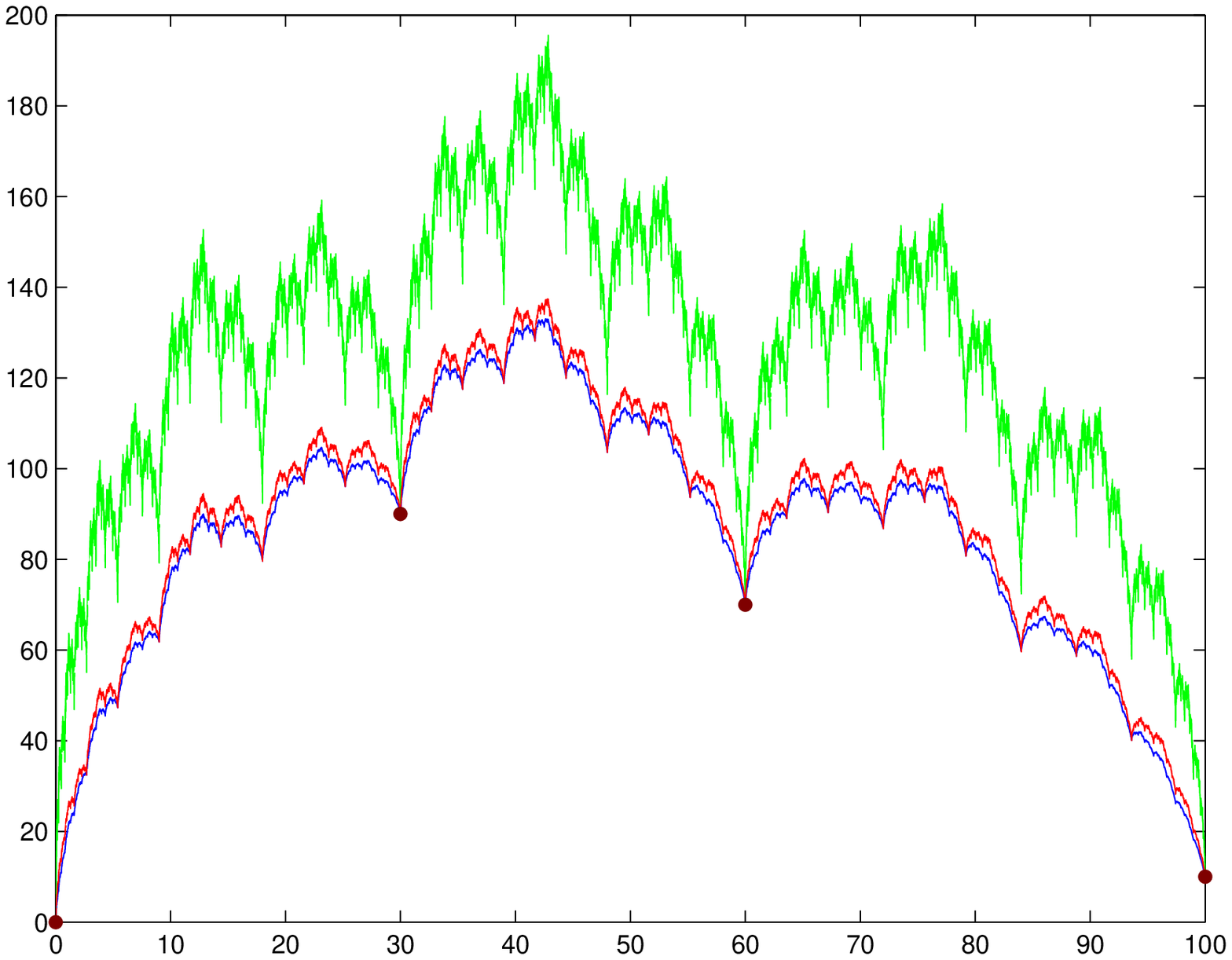, width=3cm}\label{fig:sf1}}
 \hspace{1cm}\subfigure[  $\sigma_{(r)} =
\overline{221}$]{\epsfig{file=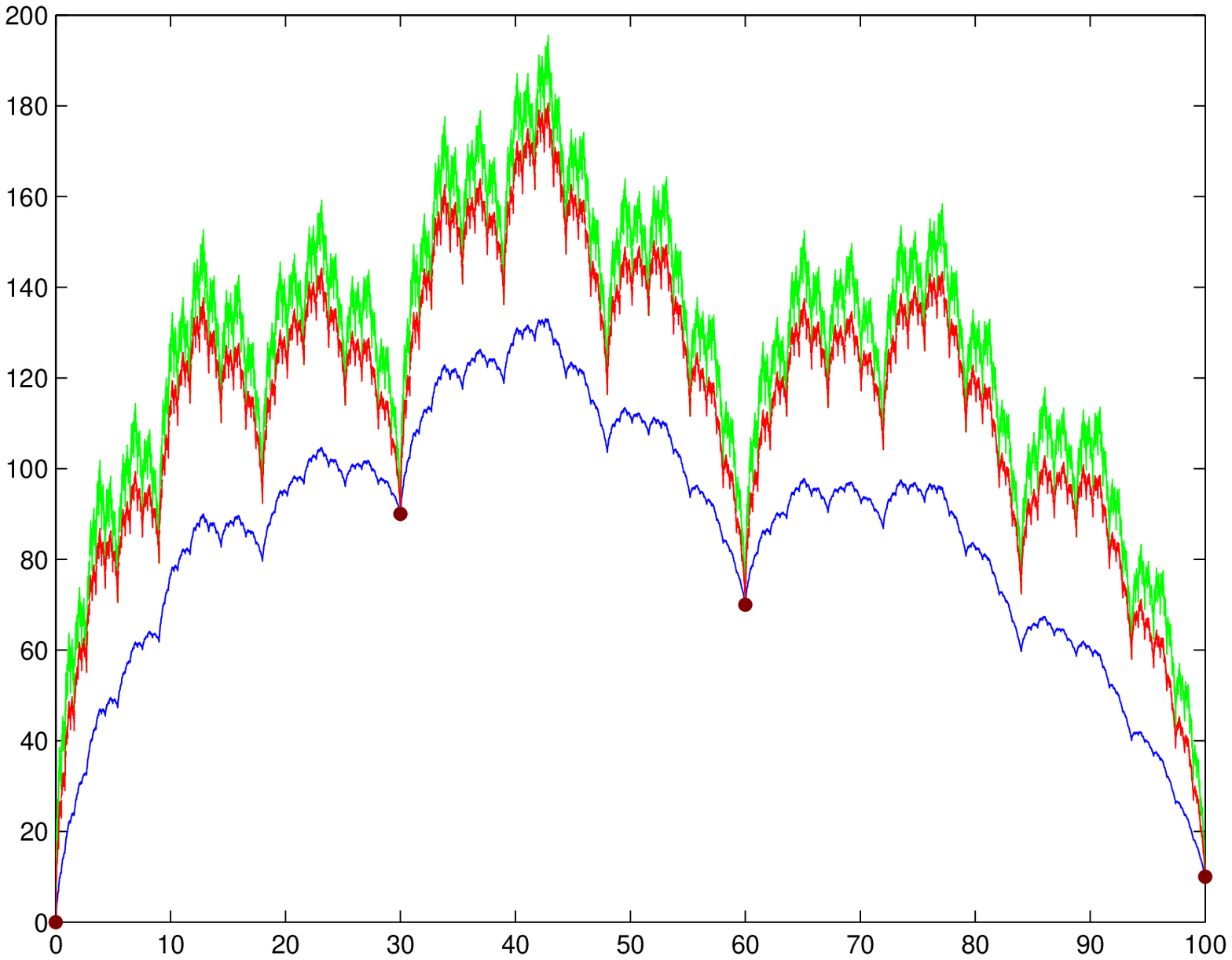, width=3cm}\label{fig:sf2}}
\hspace{1cm} \subfigure[  $\sigma_{(r)} =
\overline{12}$]{\epsfig{file=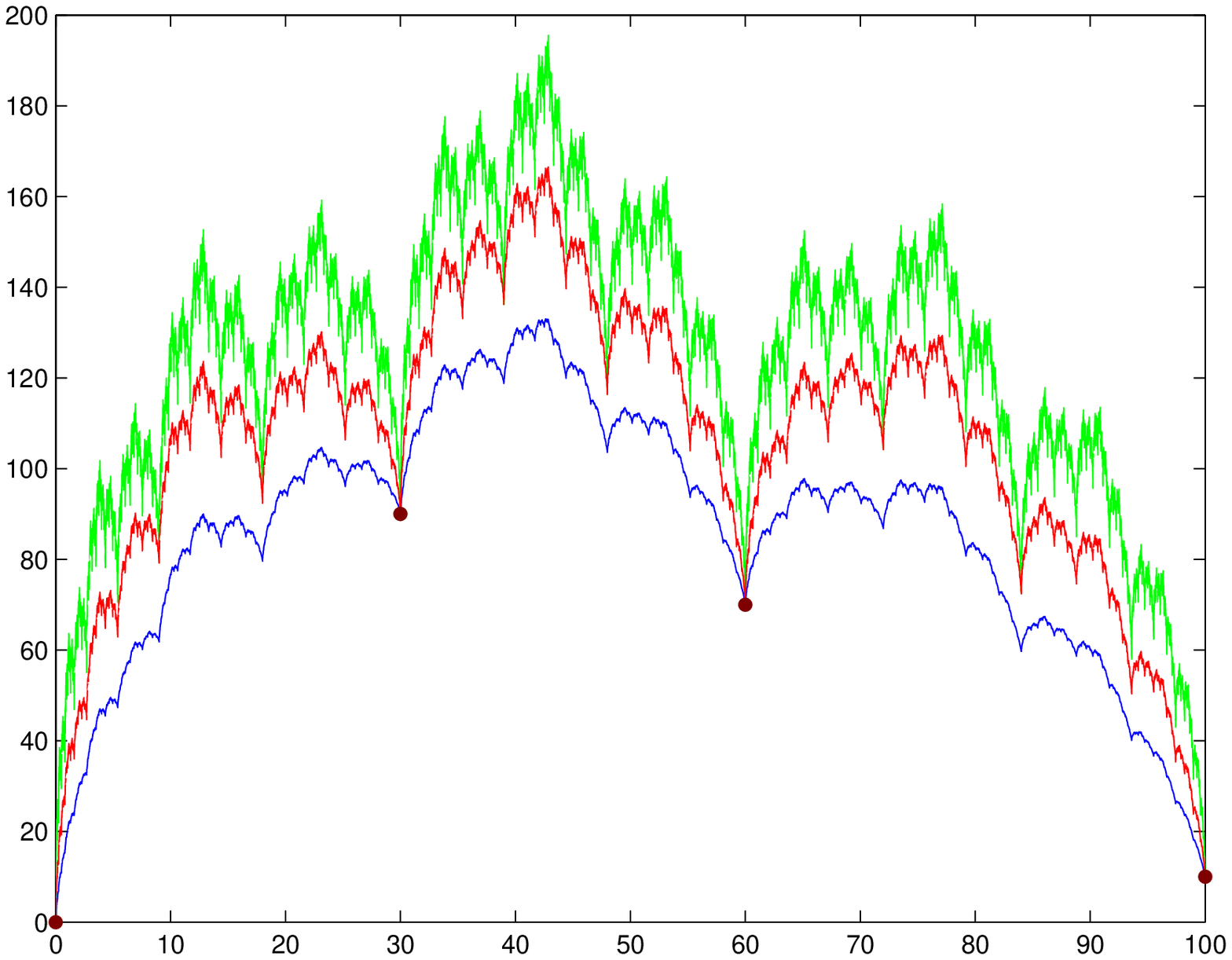, width=3cm}\label{fig:sf3}}
  \\
} {\footnotesize
\begin{flalign*}
& \mbox{Blue Curve (\textcolor{blue}{\textbf{---}}):} \ SFIF
g_{\sigma_{(b)}}
 \ \mbox{for} \  \sigma_{(b)} = \bar{1} \
\mbox{in} (a), (b) \ \mbox{and} \ (c),\\
& \mbox{Green Curve (\textcolor{green}{\textbf{---}}):} \ SFIF
g_{\sigma_{(g)}}
 \ \mbox{for} \  \sigma_{(g)} = \bar{2} \
\mbox{in} (a), (b) \ \mbox{and} \ (c),\\
& \mbox{Red Curve (\textcolor{red}{\textbf{---}}):} \ SFIF
g_{\sigma_{(r)}}
 \ \mbox{for} \  \sigma_{(r)} = \overline{112} \
\mbox{in} (a),
 \ \mbox{for} \  \sigma_{(r)} = \overline{221} \
\mbox{in} \ (b) \ \mbox{and}
 \ \mbox{for} \  \sigma_{(r)} = \overline{12} \
\mbox{in} \ (c)
\end{flalign*}}
\caption{{\small SFIFs for $\sigma_{(b)} = \bar{1}$, $\sigma_{(g)} =
\bar{2}$ and different choices of $\sigma_{(r)}$ }}
\end{figure}

In the construction of SFIF for a  $\sigma =\sigma_1\sigma_2 \ldots
\in \Lambda$, the set  $S_j = W_{\sigma_j}(S_{j-1}), \ j~=~1,
2,\ldots,$  representing the action of SIFS~\eqref{eq:sifs} at the
iteration level $j$ is computed.  The SFIF $g_{\sigma_{(b)}}$ for
$\sigma_{(b)} = \bar{1}$ (c.f. Figs.~\ref{fig:sf1}-~\ref{fig:sf3},
blue curve) is constructed by the action of IFS $\{\mathbb{R}^2 ;
\om_{n,1} , n=1,\ldots,N\} $  at every level of iteration.
Similarly, SFIF $g_{\sigma_{(g)}}$ for $\sigma_{(g)} = \bar{2}$
(c.f. Figs.~\ref{fig:sf1}-~\ref{fig:sf3}, green curve) is
constructed by the action of IFS $\{\mathbb{R}^2 ; \om_{n,2} ,
n=1,\ldots,N\} $  at every level of iteration. The SFIF
$g_{\sigma_{(r)}}$ for $\sigma_{(r)} = \overline{112}$ (c.f.
Fig.~\ref{fig:sf1}, red curve) is constructed by the action of IFS
$\{\mathbb{R}^2 ; \om_{n,1} , n=1,\ldots,N\} $ at $j^{th}$ level of
iteration if $j $ is not divisible by $3$ and by the action of IFS
$\{\mathbb{R}^2 ; \om_{n,2} , n=1,\ldots,N\} $ if $j$ is divisible
by $3$. Likewise, SFIF $g_{\sigma_{(r)}}$ for $\sigma_{(r)} =
\overline{221}$ (c.f. Fig.~\ref{fig:sf2}, red curve) is constructed
by the action of IFS $\{\mathbb{R}^2 ; \om_{n,1} , n=1,\ldots,N\} $
at $j^{th}$ level of iteration if $j$ is divisible $3$  and
otherwise by the action of IFS $\{\mathbb{R}^2 ; \om_{n,2} ,
n=1,\ldots,N\} $. Finally, SFIF $g_{\sigma_{(r)}}$ for $\sigma_{(r)}
= \overline{12}$ (c.f. Fig.~\ref{fig:sf3}, red curve) is constructed
by the action of IFS $\{\mathbb{R}^2 ; \om_{n,1} , n=1,\ldots,N\} $
at $j^{th}$ level of iteration if $j$ is not divisible by $2$ and by
the action of IFS $\{\mathbb{R}^2 ; \om_{n,2} , n=1,\ldots,N\} $ if
$j$ is divisible by $2$.

\vspace{1cm}

\begin{table}[!hbp]
\begin{center}
\begin{tabular}{|c|c|c|c|}  \hline  & i=1 & i=2 & i=3 \\ \hline
$a$    & 0.3  &0.3   &0.4   \\  \hline $b$    & 0    &30    &60
\\ \hline $e_{i,1}$  & 0.86 &-0.24 &-0.64 \\ \hline $f_{i,1}$  & 0    &90
&70    \\ \hline $e_{i,2}$  & 0.84 &-0.26 &-0.66 \\ \hline $f_{i,2}$
& 0 &90    &70    \\ \hline
\end{tabular}
\end{center}\caption{Computed Values of
$a_i,b_i,e_{i,1},f_{i,1},e_{i,2},f_{i,2}$,  $i=1,2,3$, for sample
data $S_0$
   \label{tab:Table-coef}}
\end{table}



 The SFIFs $g_{\sigma_{(b)}}
 \ \mbox{for} \  \sigma_{(b)} = \bar{1}$ and  $g_{\sigma_{(g)}}
 \ \mbox{for} \  \sigma_{(g)} = \bar{2}$   are in fact FIFs
(c.f. Figs.~\ref{fig:sf1}-~\ref{fig:sf3}, blue and green curves ),
since these are constructed with a single element of
SIFS~\eqref{eq:sifs}.

Heuristically, in terms of their fractal
dimension~\cite{barnsley86}, the graphs of SFIF $g_{\sigma_{(r)}}$
appear to fill more space in $\mathbb{R}^2$ than the graph of FIF
$g_{\sigma_{(b)}}$ and less space in $\mathbb{R}^2$ than the graph
of FIF $g_{\sigma_{(g)}}$. In fact, the fractal dimension of graphs
of FIF $g_{\sigma_{(b)}}$ and $g_{\sigma_{(g)}}$ are computed as
$1.3069$ and $1.5199$ respectively whereas the fractal dimension of
SFIF $g_{\sigma_{(r)}}$ with $\sigma_{(r)} = \overline{112}$ (c.f.
Fig.~\ref{fig:sf1}, red curve)  is  $1.3632$,  the fractal dimension
of SFIF $g_{\sigma_{(r)}}$ with $\sigma_{(r)} = \overline{221}$
(c.f. Fig.~\ref{fig:sf2}, red curve) is  $1.4572$ and the fractal
dimension of SFIF $g_{\sigma_{(r)}}$ with $\sigma_{(r)} =
\overline{12}$ (c.f. Fig.~\ref{fig:sf3}, red curve) is $1.4182$.

Further,  for   FIFs $g_{\sigma_{(b)}}$ and $g_{\sigma_{(g)}}$, the
\textit{average fractal distance} defined as $\\ d_F(f,g) =
\frac{1}{(b-a)} \left(\int\limits_a^b |f(x)-g(x)|^2 dx \right)^{1/2}
$ for the functions $f$ and $g$, continuous on a closed
interval~$[a,b]$, is $d_F(g_{\sigma_{(b)}},g_{\sigma_{(g)}}) =
0.297$.   It is observed that (i) for SFIF $g_{\sigma_{(r)}}$ with $
\sigma_{(r)} = \overline{112}$,
$d_F(g_{\sigma_{(b)}},g_{\sigma_{(r)}}) = 0.022$  while
$d_F(g_{\sigma_{(g)}},g_{\sigma_{(r)}}) = 0.276$.   So, if the data
generating function is at one third average fractal distance from
FIF $g_{\sigma_{(b)}}$, then SFIF $g_{\sigma_{(r)}}$  is a better
approximation of the data generating function, since
$g_{\sigma_{(r)}}$ is closer to $g_{\sigma_{(b)}}$ than
$g_{\sigma_{(g)}}$ (c.f. Fig.~\ref{fig:sf1}) i.e.
$d_F(g_{\sigma_{(b)}},g_{\sigma_{(r)}}) <
d_F(g_{\sigma_{(g)}},g_{\sigma_{(r)}})$. (ii) For SFIF
$g_{\sigma_{(r)}}$  with $ \sigma_{(r)} = \overline{221}$,
$d_F(g_{\sigma_{(b)}},g_{\sigma_{(r)}}) = 0.228 $ while
$d_F(g_{\sigma_{(g)}},g_{\sigma_{(r)}})  = 0.071$. So, if the data
generating function is at one third average fractal distance from
FIF $g_{\sigma_{(g)}}$, then SFIF $g_{\sigma_{(r)}}$  is a better
approximation of such data generating function, since
$g_{\sigma_{(r)}}$ is  closer to $g_{\sigma_{(g)}}$ than
$g_{\sigma_{(b)}}$ (c.f. Fig.~\ref{fig:sf2}) and (iii) for
$g_{\sigma_{(r)}}$ with $ \sigma_{(r)} = \overline{12}$,
$d_F(g_{\sigma_{(b)}},g_{\sigma_{(r)}}) = 0.138$ and
$d_F(g_{\sigma_{(g)}},g_{\sigma_{(r)}})  = 0.162$.  So, if the data
generating function lies in the middle of FIFs $g_{\sigma_{(b)}}$
and $g_{\sigma_{(g)}}$,  then SFIF $g_{\sigma_{(r)}}$ (c.f.
Fig.~\ref{fig:sf3})  is a better approximation of such data
generating function.

\section{Integral and Derivative of SFIF}\label{sec:idsfif}

In this section, for a SFIF passing through a given interpolation
data, its integral is shown to be also a SFIF, albeit for a
different interpolation data. Further, in this section, the
smoothness of SFIF is investigated in terms of its Lipschitz
exponent and it is found that, in general, a SFIF  may not be
differentiable. This, as a natural follow up, led to determining in
this section the conditions for existence of derivatives of a SFIF.

In order to study the integral of a SFIF,  a SIFS
\begin{align}\label{eq:gsifs}
\Big\{ \big\{\mathbb{R}^2;\ \om_{n,k}(x,y) = (L_n(x),G_{n,k}(x,y)) :
n = 1,\ldots,N \big\}, \ k = 1,\ldots,M \Big\},
\end{align}
associated with the data $\{ (x_i,y_i) \in \mathbb{R}^2 :i=0,
\ldots, N \}$ is considered, where $ \ L_n(x) = a_n x + b_n $ are
given by~\eqref{eq:Ln} and the functions $\ G_{n,k}(x,y)$ defined by
\begin{align}\label{eq:Gnk}
G_{n,k}(x,y) = \gm_{n,k} y + q_{n,k} (x)  ,\quad  n = 1, \ldots, N.
\end{align}
satisfy the join up conditions given by~\eqref{eq:jc}. Here,
$\gm_{n,k}$ are free parameters chosen such that $|\gm_{n,k}| < 1$
and $\gm_{n,k} \neq \gm_{n,l}$ for $k \neq l$ and $ q_{n,k} (x)$ are
affine functions. Condition~\eqref{eq:jc} ensures that there exits a
unique attractor $G_{\sigma} \in {\cal H}(\mathbb{R}^2)$ of
SIFS~\eqref{eq:gsifs}. By the arguments similar to those in the
proof of Theorem~\ref{th:SFIF}, $G_{\sigma}$ is graph of a
continuous  function $g_{\sigma}$.

The following notations~\cite{barnsley89_1} are needed in the sequel
for tidy presentation of our results:
\begin{align} \label{eq:qycap} \left. \begin{array}{ll}
 \hat{\gm}_{n,k} & = a_n \gm_{n,k}   \\
 \hat{y}_{N,k} & =  \hat{y}_0  + \frac{\sum\limits_{j=1}^N a_j \big[
\int\limits_{x_0}^{x_N} q_{j,k} (t) \ dt \big]}{1 -
\sum\limits_{j=1}^N a_j \gm_{j,k} }   \\
 \hat{y}_{n,k} & = \hat{y}_0 + \sum\limits_{j=1}^n a_j \big[ \gm_{j,k}
(\hat{y}_{N,k}-\hat{y}_0) \ + \int\limits_{x_0}^{x_N} q_{j,k} (t) \
dt
\big]   \\
 \hat{q}_{n,k}(x) & =  \hat{y}_{n-1,k} - a_n \gm_{n,k}
\hat{y}_0  + a_n \int\limits_{x_0}^x q_{n,k}(t) \ dt \end{array}
\right\}
\end{align}
 where, $\hat{y}_0$ is an arbitrary real number. To
determine an interpolation data through which the integral of SFIF
passes, let the affine functions $ q_{n,k} (x)$ in~\eqref{eq:Gnk}
satisfy :
\begin{align}\label{eq:alqrl}
\frac{\sum\limits_{j=1}^N a_j \int\limits_{x_0}^{x_N} q_{j,k}}{1 -
\sum\limits_{j=1}^N a_j \gm_{j,k}} & = \frac{\sum\limits_{j=1}^N a_j
\int\limits_{x_0}^{x_N} q_{j,l}}{1 - \sum\limits_{j=1}^N a_j
\gm_{j,l}} \neq 1 \quad \mbox{for } \ k \neq l,\ k,l = 1,\ldots,M.
\end{align}
For example, for $a_j = \frac{1}{N}$,  $\gm_{j,k} = \gm_k$ and
$q_{j,k} = (1-\gm_k)(e_j x + f_j)$  for $j=1,\ldots,N$, the
condition~\eqref{eq:alqrl} is satisfied.  Then, $\hat{y}_{i,k} =
\hat{y}_{i,l} = \hat{y}_i $ for $ i =0,\ldots,N$; $k,l = 1,\ldots,M$
and $\hat{y}_N -\hat{y}_0 \neq 1$.

The SIFS associated with the data $\{ (x_i,\hat{y}_i) \in
\mathbb{R}^2 :i=0, \ldots N \}$ is now defined as the pool of IFS
\begin{align}\label{eq:igsifs}
\Big\{ \big\{\mathbb{R}^2;\ \hat{\om}_{n,k}(x,y) =
(L_n(x),\hat{G}_{n,k}(x,y)) :\ n = 1,\ldots,N \big\}, k = 1,\ldots,M
\Big\}
\end{align}
where, the functions
\begin{align}\label{eq:G2}
\hat{G}_{n,k}(x,y) & = \hat{\gm}_{n,k} y + \hat{q}_{n,k}(x)
\end{align}
satisfy the join-up conditions $\ \hat{G}_{n,k}(x_0, \hat{y}_0) =
\hat{y}_{n-1}\ $ and $\ \hat{G}_{n,k}(x_N, \hat{y}_N)  = \hat{y}_n
$. These join-up conditions ensure that there exits a unique
attractor $\hat{G}_{\sigma} \in {\cal H}(\mathbb{R}^2)$ of
SIFS~\eqref{eq:igsifs}. The following theorem shows that the
integral of SFIF is also a SFIF albeit for interpolation data $\{
(x_i,\hat{y}_i) \in \mathbb{R}^2 :i=0, \ldots N \}$.

\begin{theorem}\label{th:intg} For the
interpolation data $\{ (x_i,y_i) \in \mathbb{R}^2 :i=0, \ldots N
\}$, let $g_{\sigma}$ be SFIF corresponding to SIFS~\eqref{eq:gsifs}
for $\sigma \in \Lambda$. Then, the integral
\begin{align}\label{eq:intg}
\hat{g}_{\sigma}(x) =  \hat{y}_0 + \int\limits_{x_0}^x g_{\sigma}(t)
\ dt
\end{align}
is SFIF associated with SIFS~\eqref{eq:igsifs} for the interpolation
data $\{ (x_i, \hat{y}_i ) : i =0,\ldots,N \}$.
\end{theorem}

\begin{proof}
Using~\eqref{eq:intg} and~\eqref{eq:rb}, it is observed that,
\begin{align}\label{eq:intg1}
\hat{g}_{\sigma}(L_{i_k} \circ \ldots \circ L_{i_1}(x))
 & =
\hat{g}_{\sigma}(L_{i_k} \circ \ldots \circ L_{i_1}(x_0))  +
\left(\prod\limits_{j=1}^k a_{i_j} \gm_{i_j,\sigma_j} \right)
\Big( \hat{g}_{\sigma}(x) -\hat{y}_0   \Big) \nno \\
& \quad \mbox{}+ \sum\limits_{p=1}^k \left(\prod\limits_{j=p+1}^k
a_{i_j} \gm_{i_j,\sigma_j} \right) a_{i_p} \int\limits_{L_{i_{p-1}}
\circ \ldots \circ L_{i_1}(x_0)}^{L_{i_{p-1}} \circ \ldots \circ
L_{i_1}(x)} q_{i_p,\sigma_p} (t).
\end{align}

Also, by~\eqref{eq:intg},
\begin{align*}
\hat{g}_{\sigma}(L_{i_k} \circ \ldots \circ L_{i_1}(x_0)) & =
\hat{y}_0 +  \sum\limits_{p=1}^k \bigg(\prod\limits_{j=p+1}^k
a_{i_j} \gm_{i_j,\sigma_j} \bigg) \bigg\{\sum\limits_{l=1}^{i_p-1}
a_l   \bigg[ \gm_{l,\sigma_p} (\hat{y}_N-\hat{y}_0) \\ & \quad
\mbox{} + \int\limits_{x_0}^{x_N} q_{l,\sigma_p}(t)\ dt \bigg]  +
a_{i_p} \int\limits_{x_0}^{L_{i_{p-1}} \circ \ldots \circ
L_{i_1}(x_0)} q_{i_p,\sigma_p}(t) \bigg\}.
\end{align*}

The above identity and~\eqref{eq:qycap} give
\begin{align}\label{eq:gL0}
\lefteqn{\hat{g}_{\sigma}(L_{i_k} \circ \ldots \circ L_{i_1}(x_0))}
\nno \\  & = \bigg(\prod\limits_{j=1}^k a_{i_j} \gm_{i_j,\sigma_j}
\bigg) \hat{y}_0  \nno \\ & \quad \mbox{}  + \sum\limits_{p=1}^k
\bigg(\prod\limits_{j=p+1}^k a_{i_j} \gm_{i_j,\sigma_j} \bigg)
\bigg\{  \hat{y}_{i_p - 1}
    -a_{i_p} \gm_{i_p,\sigma_p}\ \hat{y}_0  + a_{i_p} \int\limits_{x_0}^{L_{i_{p-1}} \circ \ldots \circ
L_{i_1}(x_0)} q_{i_p,\sigma_p}(t) \bigg\}.
\end{align}

Now, substituting the value of $\hat{g}_{\sigma}(L_{i_k} \circ
\ldots \circ L_{i_1}(x_0))$ from~\eqref{eq:gL0} in~\eqref{eq:intg1},
it follows that
\begin{align*}
\hat{g}_{\sigma}(L_{i_k} \circ \ldots \circ L_{i_1}(x))
 & = \hat{G}_{i_k,\sigma_{k}} \bigg( L_{i_{k-1}} \circ
\ldots \circ L_{i_1} (x), \hat{G}_{i_{k-1},\sigma_{k-1}} \Big( .,
\ldots  \nno \\ & \quad \quad \quad \quad  \hat{G}_{i_2,\sigma_2}
\big(L_{i_1}(x) , \hat{G}_{i_1,\sigma_1}(x,\hat{g}_{\sigma}(x))\big)
\ldots \Big)\bigg).
\end{align*}
Thus,  $\hat{g}_{\sigma}$ is SFIF associated with
SIFS~\eqref{eq:igsifs}. \qed \end{proof}

\begin{remark}\label{rem:psfifint}
In case of $\kp$-SFIF $ g_{\sigma,\kp}$ (c.f.
Remark~\ref{rem:kpsfif}),  using the lines of proof of
Theorem~\ref{th:intg}, it follows that the integral of $\kp$-SFIF is
not a $\kp$-SFIF but integral of $g_{\sigma,\kp} + \xi_{\sigma}$,
where $\xi_{\sigma}$ is defined by
\begin{align}\label{eq:xi}
 \xi_{\sigma}(x) = \gm_{i_n,\sigma_n}
\left[ \xi_{\sigma} (L_n^{-1} (x)) - \kp ( 1- L_n^{-1}(x))\right] +
\kp (1-x)
\end{align}
for $ x \in I_n$,  is a $\kp$-SFIF for the interpolation data $\{
(x_i,\hat{y}_i) \in \mathbb{R}^2 :i=0, \ldots N \}$, provided
\begin{align*}
\frac{\sum\limits_{j=1}^N a_j [\int\limits_{x_0}^{x_N} (q_{j,k}(t) -
\kp \gm_{j,k}  (1-t)) dt]}{1 - \sum\limits_{j=1}^N a_j \gm_{j,k}} &
= \frac{\sum\limits_{j=1}^N a_j [\int\limits_{x_0}^{x_N} (q_{j,l}(t)
- \kp \gm_{j,l}  (1-t)) dt]}{1 - \sum\limits_{j=1}^N a_j \gm_{j,l}}
\neq 1 \quad  \mbox{holds}.
\end{align*}
Here,  $\hat{y}_{j,k}$ , $\hat{q}_{j,k}(x),\ j=1,2,\ldots, N$,
in~\eqref{eq:qycap} are given by
\begin{align*}
(1-\kp) \hat{y}_{N,k} & = \kp (x_0-x_N) + (1-\kp) \hat{y}_0 \nno \\
&  \mbox{} + \frac{\sum\limits_{n=1}^N a_n \bigg\{
\int\limits_{x_0}^{x_N} q_{n,k} (t) \ dt + \kp
\int\limits_{x_0}^{x_N} (1 - L_n(t)) \ dt  - \kp\ \gm_{n,k}
\int\limits_{x_0}^{x_N} (1-t) \ dt \bigg\}}{1 - \sum\limits_{n=1}^N
a_n \gm_{n,k} }
\end{align*}
\begin{align*}
 (1-\kp) \hat{y}_{j,k} & = \kp (x_0-x_j) + (1-\kp)
\hat{y}_0 + \sum\limits_{n=1}^j a_n \bigg[ \gm_{n,k} \bigg\{ (1-\kp)
(\hat{y}_{N,k}-\hat{y}_0)\nno \\ & \quad \mbox{} + \kp
\int\limits_{x_0}^{x_N} t \ dt \bigg\}
 +  \int\limits_{x_0}^{x_N} \{\kp (1 - L_n(t))  +  q_{n,k} (t)\} \ dt
 \bigg]
\end{align*}
and
\begin{align}\label{eq:qycapksfif}
\hat{q}_{j,k}(x) & = \kp x_{j-1} + (1-\kp) \hat{y}_{j-1,k} - a_j
\gm_{j,k} \bigg[\kp x_0 + (1-\kp) \hat{y}_0 + \kp
\int\limits_{x_0}^x (1-t) \ dt \bigg] \nno
\\ & \quad \mbox{} + a_j \kp \int\limits_{x_0}^x (1 - L_j(t))\ dt +
a_j \int\limits_{x_0}^x q_{j,k}(t) \ dt .
\end{align}
\end{remark}

For investigating the smoothness of a SFIF, the following notations
are needed:
\begin{align}\label{eq:laC1}
\la &= \min \{\la_{n,k} : n=1,2,\ldots,N , \ k =1,2,\ldots M
\},  \mbox{where} \ \la_{n,k} \ \mbox{are real numbers} \nno \\
&  \hspace{8cm} \mbox{
satisfying} \ 0 < \la_{n,k} \leq 1 \nno \\
C_1 &= \max \{\frac{|\gm_{n,k}|}{|I_n|^{\la}} : n=1,2,\ldots,N,\
k=1,2,\ldots,M \},  \mbox{where} \ \gm_{n,k} \ \mbox{are real
numbers}
 \nno \\ &  \hspace{8cm}  \mbox{satisfying} \ | \gm_{n,k}| \leq 1 .  \nno \\ &
 \mbox{Modulus of continuity of }\ g_{\sigma}(x)\ \mbox{as} \
\omega(g_{\sigma},t) = \max\limits_{ |h| \leq t} \max\limits_x
|g_{\sigma}(x + h) -g_{\sigma}(x)|
\end{align}

 The smoothness of a SFIF  in terms of its Lipscitz
exponent is given by the following theorem:

\begin{theorem}\label{th:smsfif}
Let $g_{\sigma}$ be a SFIF corresponding to SIFS~\eqref{eq:gsifs}
with $q_{n,k} \in \mbox{Lip} \ \la_{n,k}$, $ 0 < \la_{n,k} \leq 1$.
  Then,
\begin{enumerate}[(i)]
 \item  for $C_1 < 1$,  $\ g_{\sigma}
\in Lip \ \la$
\item   for $C_1 = 1$,  $\ \omega(g_{\sigma},t)=O(|t|^{\la} \log |t|)$
\item for $C_1 > 1$,   $\ g_{\sigma} \in Lip \ \bar{\la}$,
\end{enumerate}
where, $\bar{\la} \leq \max\limits_{\substack{n=1,\ldots,N \\
k=1,\ldots ,M}} \left(\frac{\log \gm_{n,k}}{\log a_n}\right) $  and
$C_1$, $\la$ are given by~\eqref{eq:laC1}.
\end{theorem}

\begin{proof}
The method of proof is similar to that in~\cite{gang96}, wherein
$\gm_n$ is replaced by $\gm_{n,\sigma_n}$. \qed \end{proof}

\begin{remark}
It follows from Theorem~\ref{th:smsfif} that $g_{\sigma} \in Lip \
\overline{\la_k}$ for $C_1 \geq C_{1,k} > 1$, where $C_{1,k} = \max
\{\frac{|\gm_{n,k}|}{|I_n|^{\la_k}} : n=1,2,\ldots,N\}$ and
$\overline{\la_k} \leq \max \{\frac{\log |\gm_{n,k}|}{\log |I_n|}  :
n=1,2,\ldots,N\}$.
\end{remark}

\begin{remark}
In case of $\kp$-SFIF  $ g_{\sigma,\kp}$(c.f.
Remark~\ref{rem:kpsfif}), the smoothness result analogous to
Theorem~\ref{th:smsfif} can be obtained as follows: (i)
$g_{\sigma,\kp} \in Lip \ \la$ for $C_1 < 1$, (ii) $\
\omega(g_{\sigma, \kp},t)=O(|t|^{\la} \log |t|)$ for $C_1 = 1$ and
(iii) $g_{\sigma,\kp} \in Lip \ \bar{\la}$  for $C_1
> 1$,     $ \bar{\la} \leq \max \{(\log
\gm_{n,k})/(\log a_n) : n=1,2,\ldots,N, \ k=1,2,\ldots ,M\}$.
\end{remark}

In general, a SFIF belonging to certain Lipschitz class, need not be
differentiable. This, as a natural follow up, leads to
identification of conditions for the existence of derivative of a
SFIF in the following proposition:

\begin{proposition}\label{prop:dfsfif}
For the interpolation data $\{ (x_i,y_i) \in \mathbb{R}^2 :i=0,1,
\ldots N \}$, let $g_{\sigma}$ be a SFIF corresponding to
SIFS~\eqref{eq:gsifs} for $\sigma \in \Lambda$. Then,
$\hat{g}_{\sigma}'$  exists and $ \hat{g}_{\sigma}'(x) =
g_{\sigma}(x) $ if and only if  $\ \hat{g}_{\sigma}$ is a SFIF
associated with SIFS~\eqref{eq:igsifs} for the interpolation data
$\{ (x_i, \hat{y}_i ) : i =0,1,\ldots,N \}$, provided
$\hat{\gm}_{j,k} = a_j \gm_{j,k} $ and $
\frac{d}{dx}(\hat{q}_{j,k}(x)) = a_j q_{j,k} (x)$ hold.
\end{proposition}
\begin{proof}
If $\ \hat{g}_{\sigma}'(x) = g_{\sigma}(x)   $, then $\
\hat{g}_{\sigma}(x) =   \hat{y}_0 + \int\limits_{x_0}^x
g_{\sigma}(t)\ dt  $,  so that  `` if '' part follows from
Theorem~\ref{th:intg}. Conversely, suppose $\ \hat{g}_{\sigma}\ $ is
a SFIF associated with SIFS~\eqref{eq:igsifs} for the interpolation
data $\{ (x_i,\hat{y}_i ) : i =0,\ldots,N \}$. Then,
\begin{align}\label{eq:qd1}
\hat{g}_{\sigma}(L_{i_k} \circ \ldots \circ L_{i_1}(x)) & =
\left(\prod\limits_{j=1}^k  \hat{\gm}_{i_j,\sigma_j} \right)
\hat{g}_{\sigma}(x) \nno \\ & \quad \mbox{}+ \sum\limits_{p=1}^k
\left(\prod\limits_{j=p+1}^k \hat{\gm}_{i_j,\sigma_j} \right)
\hat{q}_{i_p,\sigma_p} (L_{i_{p-1}} \circ \ldots \circ L_{i_1} (x)).
\end{align}

Since, $ \frac{d}{dx}(\hat{q}_{j,k}(x)) = a_j q_{j,k} (x)$
\begin{align}\label{eq:qjkcap}
\hat{q}_{j,k}(x) &= \hat{q}_{j,k}(x_0) + a_j \int\limits_{x_0}^{x}
q_{j,k}(t) dt  = \hat{y}_{j-1} - a_j \gm_{j,k} \hat{y}_0 + a_j
\int\limits_{x_0}^{x} q_{j,k}(t) dt .
\end{align}

Substituting~\eqref{eq:qjkcap} and $\hat{\gm}_{j,k} = a_j \gm_{j,k}$
in~\eqref{eq:qd1},
\begin{align}\label{eq:convg}
\hat{g}_{\sigma}(L_{i_k} \circ \ldots \circ L_{i_1}(x)) & =
\Big(\prod\limits_{j=1}^k  a_{i_j} \gm_{i_j,\sigma_j} \Big)
\hat{g}_{\sigma}(x)  + \sum\limits_{p=1}^k
\Big(\prod\limits_{j=p+1}^k a_{i_j}\gm_{i_j,\sigma_j} \Big) \times
\nno
\\ & \quad \mbox{} \times \Big[\hat{y}_{i_p-1} - a_{i_p} \gm_{i_p,\sigma_p}
\hat{y}_0 + a_{i_p} \int\limits_{x_0}^{L_{i_{p-1}} \circ \ldots
\circ L_{i_1} (x)} q_{i_p,\sigma_p}(t) dt \Big] .
\end{align}
For a fixed $\sigma \in \Lambda $, it is easily seen that the
Read-Bajraktarevic operator $\hat{T}$ defined by
\begin{align}\label{eq:rb2}
\hat{T}(\sigma, g)(x) & =  \lim\limits_{k \rightarrow \infty}
\bigg\{\Big(\prod\limits_{j=1}^k a_{i_j} \gm_{i_j,\sigma_j} \Big)
g(L_{i_1}^{-1} \circ \ldots \circ L_{i_k}^{-1}(x))  +
\sum\limits_{p=1}^k \Big(\prod\limits_{j=p+1}^k
a_{i_j}\gm_{i_j,\sigma_j} \Big) \times \nno
\\ & \quad \mbox{} \times \Big[\hat{y}_{i_p-1} - a_{i_p} \gm_{i_p,\sigma_p}
\hat{y}_0 + a_{i_p} \int\limits_{x_0}^{L_{i_p}^{-1} \circ \ldots
\circ L_{i_k}^{-1} (x)} q_{i_p,\sigma_p}(t) dt \Big]\bigg\}
\end{align}
 is a
contraction map on $  {\cal \hat{G}} = \{ f : I \rightarrow
\mathbb{R} \ \ \mbox{such that} \ f \ \mbox{ is continuous},\ f(x_0)
= \hat{y}_0 \ \mbox{and} \\ f(x_N) = \hat{y}_N\} $.
By~\eqref{eq:convg}, the function $\hat{g}_{\sigma}$ is a fixed
point of~$\hat{T}$. Also, by Theorem~\ref{th:intg}, the function
$h(x) =\hat{y}_0 + \int\limits_{x_0}^x g_{\sigma}(t)\ dt \ $ is  a
SFIF associated with SIFS~\eqref{eq:igsifs}
satisfying~\eqref{eq:convg}. Consequently, $h$ also is a fixed point
of~$\hat{T}$. Hence, by uniqueness of fixed point of
Read-Bajraktarevic operator $\hat{T}$,  $\ \hat{g}_{\sigma}(x) =
\hat{y}_0 + \int\limits_{x_0}^x g_{\sigma}(t)\ dt \ $ which implies
that $\hat{g}_{\sigma}'$  exists and $ \hat{g}_{\sigma}'(x) =
g_{\sigma}(x) $, since $g_{\sigma}$ being a SFIF corresponding to
SIFS~\eqref{eq:gsifs}, is a continuous function. \qed \end{proof}

For the investigation of $n^{th}$ derivative of SFIF,  denote
\begin{align}\label{eq:Gderi}
G_{i,k,j}(x,y)  = \gm_{i,k,j}\ y  + q_{i,k,j} (x)
\end{align}
where,  $ G_{i,k,0}(x,y) = G_{i,k}(x,y)$, $\ q_{i,k,0} (x) = q_{i,k}
(x)$, $\ \gm_{i,k,0} = \gm_{i,k} \ $  and  $\ G_{i-1,k,j}(x_N,
y_{N,k,j}) = G_{i,k,j}(x_0, y_{0,k,j})  $,  $ \ i=1,\ldots,N ,\
k=1,\ldots,M \ \mbox{and} \ j = 0, 1,\ldots, n$. To determine
interpolation data through which derivatives of SFIF passes, let the
affine maps $q_{i,k,j} (x)$ in~\eqref{eq:Gderi} satisfy :
\begin{align}\label{eq:qderi}
 \frac{\sum\limits_{p=1}^N a_p \int\limits_{x_0}^{x_N} q_{p,k,j}}{1
- \sum\limits_{p=1}^N a_p \gm_{p,k,j}}  = \frac{\sum\limits_{p=1}^N
a_p \int\limits_{x_0}^{x_N} q_{p,l,j}}{1 - \sum\limits_{p=1}^N a_p
\gm_{p,l,j}} \neq 1,
\end{align}
where $\hat{y}_{0,j}, j = 0, 1,\ldots, n$ are arbitrary real
numbers. For example, for $a_i = \frac{1}{N}$, $\gm_{i,k,j} =
\gm_{k,j}$ and $q_{i,k,j} = (1-\gm_k)\bar{q}_{i,j}(x)$, where
$\bar{q}_{i,j}(x)$ are polynomials of degree $n-j$ for
$i=1,\ldots,N$, the condition~\eqref{eq:qderi} is satisfied. Then,
$\hat{y}_{i,k,j} = \hat{y}_{i,l,j} = \hat{y}_{i,j} $ for $ i
=1,\ldots,N$,  $k,l = 1,\ldots,M$ and $j = 0, 1,\ldots, n$.

The SIFS associated with the interpolation data $\{(x_i,y_{i,j}) : i
= 0,1,\ldots,N\}$,   $  j = 0,1,\ldots, n$, is now defined as
\begin{align}\label{eq:sifsderi}
\Big\{ \big\{\mathbb{R}^2;\ \om_{i,k,j}(x,y) =
(L_i(x),G_{i,k,j}(x,y)) : i = 1,\ldots,N\big\},\ k =1,\ldots,M
\Big\}.
\end{align}
It is observed that SIFS~\eqref{eq:sifsderi} reduces to
SIFS~\eqref{eq:gsifs} if $j=0$. The following theorem gives the
existence of derivatives of a SFIF.

\begin{theorem}\label{th:ndsfif}
Let  the functions $G_{i,k,j}(x,y)$ defined in~\eqref{eq:Gderi} be
such that, for some integer $ n \geq~0$,  $|\gm_{i,k}| < a_i^n$, $ \
q_{i,k} \in C^n[x_0,x_N],\ i=1,2,\ldots,N, \ k =1,2,\ldots,M $ and
$g_{\sigma}$ be a SFIF corresponding to SIFS~\eqref{eq:sifsderi} for
$j=0$  and $\sigma \in \Lambda$. Then, for $  j = 1,2,\ldots, n$,
$g_{\sigma}^{(j)} $ exists and  is a SFIF associated with
SIFS~\eqref{eq:sifsderi} for the interpolation data $\{(x_i,y_{i,j})
: i = 0,1,\ldots,N\}$, provided $ \gm_{i,k,j} =
\frac{\gm_{i,k}}{a_i^j} \ $  and $\ q_{i,k,j}(x) =
\frac{q_{i,k,j-1}^{(1)} (x) }{a_i}$.
\end{theorem}

\begin{proof}
The equation $ G_{1,k,j}\left(x_0, y_{0,j}\right) = y_{0,j}\ $ gives
$ \  y_{0,j}  = \frac{\gm_{1,k}}{a_1^j}   y_{0,j} +
 \frac{q_{1,k}^{(j)}(x_0)}{a_1^j}\ $  which implies  $ \\
y_{0,j}  = \frac{q_{1,k}^{(j)}(x_0)}{ (a_1^j - \gm_{1,k})} $.
Similarly, $\ G_{N,k,j}\left(x_N, y_{N,j}\right) = y_{N,j} \ $ gives
$ y_{N,j} = \frac{q_{N,k}^{(j)}(x_N)}{ (a_N^j-\gm_{N,k})} $. By
Proposition~\ref{prop:dfsfif}, it now follows that, for $  j =
1,2,\ldots, n$, $ g_{\sigma}^{(j)} ( x ) $  is the  SFIF associated
with SIFS $\Big\{ \big\{\mathbb{R}^2;\ \om_{i,k,j}(x,y) =
(L_i(x),G_{i,k,j}(x,y)) : i = 1,2,\ldots,N\big\},\ k =1,2,\ldots,M
\Big\}$. \qed
\end{proof}

\begin{remark}
In the case of $\kp$-SFIF, Remark~\ref{rem:psfifint} and
Theorem~\ref{th:ndsfif} suggest that the function $
\tilde{g}_{\sigma,\kp}^{(j)} = g_{\sigma,\kp}^{(j)} ( x ) -
\sum\limits_{p=1}^j \xi_{\sigma}^{(p-1)}(x)$, with $\xi_{\sigma}$
given by~\eqref{eq:xi},  is $\kp$-SFIF associated with the SIFS
\begin{align*} \Big\{ \big\{\mathbb{R}^2;\ \om_{i,k,j} : i = 1,2,\ldots,N\big\},\
k =1,2,\ldots,M \Big\},     j = 0,1,\ldots, n. \end{align*}
Using~\eqref{eq:xi} and~\eqref{eq:qycapksfif}, it is easily seen
that  $\ y_{0,1} = \frac{q_{1,k}^{(1)} (x_0)}{(1 - \kp)
(a_1-\gm_{1,k})} - \frac{\kp}{1 - \kp} \ $  and  $\\ y_{N,k,1} =
\frac{q_{N,k}^{(1)}(x_N)}{(1 - \kp) (a_N-\gm_{N,k})} - \frac{\kp}{1
- \kp}$. Also, $\ y_{0,j}  = \frac{q_{1,k}^{(j)} (x_0)}{(1 - \kp)
(a_1^j-\gm_{1,k})} \quad \mbox{and} \ y_{N,k,j} =
\frac{q_{N,k}^{(j)}(x_N)}{(1 - \kp) (a_N^j-\gm_{N,k})}$, \\$k
=1,2,\ldots,M$,  for $j > 1 $.
\end{remark}

\section{Conclusions}
In the present work, the notion of Super Fractal Interpolation
Function (SFIF) is introduced for finer simulation of the objects of
the nature or outcomes of scientific experiments that reveal one or
more structures embedded in to another. Since, in the construction
of SFIF, at each level of iteration, an IFS can be chosen from a
pool of several IFS, the desired randomness and variability can be
implemented in fractal interpolation of the given data.  Thus, SFIF
may be used as a tool for better geometrical modeling of objects
found in nature and results of certain scientific experiments. Also,
an expository description of investigations on the integral, the
smoothness and determination of conditions for existence of
derivatives  of a SFIF is given in the present work. It is proved
that, for a SFIF passing through a given interpolation data, its
integral is also a SFIF, albeit for a different interpolation data.
The smoothness of a SFIF is given in terms of its Lipschitz
exponent. A SFIF $g_{\sigma}$, for $C_1 \neq 1$,  belongs to a
Lipschitz class and, for $C_1 = 1 $, $\
\omega(g_{\sigma},t)=O(|t|^{\la} \log |t|)$. It is seen that the
smoothness of SFIF  depend on free variables $\gm_{n,k}$ as well as
on the smoothness of affine functions $q_{n,k}(x)$ occurring in its
definition. Further, sufficient conditions for existence of
derivatives of a SFIF are derived in the present paper. Our results
on SFIF found here are likely to have wide applications in areas
like pattern-forming alloy solidification in chemistry, blood vessel
patterns in biology, signal processing, fragmentation of thin plates
in engineering, stock markets in finance, wherein significant
randomness and variability is observed in simulation of various
processes.

\section*{Acknowledgement}
The second author thanks CSIR for Research Grant No:
9/92(417)/2005-EMR-I for the present work.

\bibliographystyle{plain}

\end{document}